\documentclass{amsart}
\usepackage{amsmath, amsthm, latexsym, amssymb, amsfonts}
\usepackage[pdftex]{graphicx,color} 

\newenvironment{pf}{\proof[\proofname]}{\endproof}
\theoremstyle{plain}
\newtheorem{Th}{Theorem}[section]

\newtheorem{Cor}[Th]{Corollary}
\newtheorem{Prop}[Th]{Proposition}
\newtheorem{Lemma}[Th]{Lemma}
\numberwithin{equation}{section}
\numberwithin{figure}{section}
\theoremstyle{definition}
\newtheorem{Rem}[Th]{Remark}

\newtheorem{Ex}[Th]{Example}
\newtheorem{Def}[Th]{Definition}

\newcommand{\cal}[1]{\mathcal{#1}}

\newcommand{\N}{\mathbb N}
\newcommand{\Z}{\mathbb Z}
\newcommand{\R}{\mathbb R}
\newcommand{\F}{\mathbb F}

\newcommand{\cA}{\cal A}
\newcommand{\cB}{\cal B}
\newcommand{\cC}{\cal C}

\newcommand{\cL}{\cal L}

\newcommand{\K}{{\mathbb K}}

\newcommand{\la}{\langle}
\newcommand{\ra}{\rangle}

\newcommand{\Img}{\operatorname{Im}}

\newcommand{\spn}{\operatorname{span}}

\newcommand{\Hilb}{\operatorname{Hilb}}
\newcommand{\GL}{\operatorname{GL}}
\newcommand{\Fan}{\operatorname{Fan}}
\newcommand{\Ker}{\operatorname{Ker}}

\newcommand{\rl}[1]{Lemma~\ref{L:#1}}
\newcommand{\rp}[1]{Proposition~\ref{P:#1}}
\newcommand{\rr}[1]{Remark~\ref{R:#1}}

\newcommand{\re}[1]{(\ref{e:#1})}

\newcommand{\rt}[1] {Theorem~\ref{T:#1}}
\newcommand{\rd}[1]{Definition~\ref{D:#1}}
\newcommand{\rf}[1]{Figure~\ref{F:#1}}

\begin{document}


\title[On dual toric complete intersection codes]{On dual toric complete intersection codes}
\author[Pinar Celebi Demirarslan]{Pinar Celebi Demirarslan}
\address[Pinar Celebi Demirarslan]{Bartin University\\ Bartin, Turkey}
\email{celebi\_pin@hotmail.com}
\author[Ivan Soprunov]{Ivan Soprunov}
\address[Ivan Soprunov]{Department of Mathematics\\ Cleveland State University\\ Cleveland, OH USA}
\email{i.soprunov@csuohio.edu}
\thanks{The second author is partially supported by NSA Grant H98230-13-1-0279}


\begin{abstract} In this paper we study duality for evaluation codes on intersections of $d$
hypersurfaces with given $d$-dimensional Newton polytopes, so called toric complete intersection codes.
In particular, we give a condition for such a code to be quasi-self-dual.
In the case of $d=2$ it reduces to a combinatorial condition on the Newton polygons.
This allows us to give an explicit construction of dual and quasi-self-dual toric complete intersection codes.
 We provide a list of examples over $\F_{16}$.

\end{abstract}

\maketitle


\section{Introduction}

In this paper we consider a class of evaluation codes called toric complete intersection codes.
They were  introduced in \cite{So1} and are a natural generalization of evaluation codes
on complete intersections in the projective space, previously studied by Duursma, Renter\'ia,  and Tapia-Recillas
 \cite{DRT}; Gold, Little, and Schenck \cite{GLS}; and  Ballico and Fontanari \cite{BF}.
 
 A toric complete intersection code $\cC_{S,A}$ is constructed by evaluating $d$-variate Laurent polynomials supported in
 a given lattice polytope $A$ at the set $S$ of common zeroes of $d$ Laurent polynomials with given Newton
 polytopes $P_1, \dots, P_d$. In \cite{So1}, the second author proved general bounds for the minimum distance of such codes
 in terms of $A$ and the $P_i$.  The goal of this paper is to study duality for toric complete intersection codes. 
 In particular, we give conditions on $A$ and $P_1,\dots, P_d$ when the code $\cC_{S,A}$ is quasi-self dual,
 see \rt{y-dual toric}. 
 
 When $d=2$ we give a combinatorial formula for the dimension of $\cC_{S,A}$, thus
 reducing the above mentioned conditions to purely combinatorial ones (see \rt{geom criterion}).
 We show how restrictive this condition is when the polytopes $P_i$ are similar. In fact, in this case
  a quasi-self dual $\cC_{S,A}$ exists if and only if the $P_i$ are $\GL(2,\Z)$-equivalent to an integer multiple 
  of one of 16 polygons as in \rp{unmixed}. On the other hand, \rt{summands} provides a much less restrictive
  framework for constructing the polytopes $A$ and $P_1, P_2$ which produce quasi-self dual codes $\cC_{S,A}$.
 
 The paper concludes with an algorithm for finding dual and  quasi-self dual toric complete intersection codes,
and provides with a list of examples over the finite field of 16 elements. 

\section{Preliminaries}

\subsection{Dual codes}
To set our notation we start with basic definitions from coding theory. 
Throughout the paper, $\F_q$
denotes a finite field of $q$ elements and $\F_q^*$ its multiplicative group of non-zero elements.
A subspace $\cC$ of $ \mathbb{F}_{q}^{n}$ is called a \textit{linear code}, and its elements 
${c}=(c_{1},\dotsc,c_{n})$ are called \textit{codewords}.  The number $n$ is called the {\it block-length} of $\cC$.
The {\it weight} of $c$ in $\cC$ is the number of non-zero entries in $c$.
The {\it distance} between two codewords $a$ and $b$ in $\mathcal{C}$ is the weight of $a-b\in\cC$.
The minimum distance between distinct codewords in $\cC$ is the same as the minimum weight of  non-zero codewords in $\cC$
and will be denoted by $d(\cC)$. The block-length $n$, the dimension $k=\dim(\cC)$, and the minimum
distance $d=d(\cC)$ are the parameters of $\cC$. A code with parameters $n$, $k$, and $d$ is referred to as an $[n,k,d]_{q}$-code.
The parameters of any $[n,k,d]_{q}$-code satisfy the Singleton bound: $d\leq n-k+1$. Codes that meet the Singleton bound are
called {\it maximum distance separable} (MDS) codes.

A  \textit{generator matrix} of a linear code $\cC$  is a matrix whose rows form a basis for $\cC$.
Two linear codes $\cC$ and  $\cC'$ are called \textit{equivalent} if a generator matrix of $\cC'$ is obtained by scaling and permuting 
the columns of a generator matrix of $\cC$. 
In particular, for any $x \in (\F_q^*)^n$ the code
\begin{displaymath}
\cC' = x\,\cC = \{ (x_1u_1, \dots, x_nu_n) \, | \, u\in \cC \}
\end{displaymath}  
is equivalent to $\cC$. Equivalent codes have the same parameters.

Let  $(u \cdot v)$ be the standard dot product on $\F_q^n$. Then we can define the {\it dual code} by
\begin{displaymath}
\mathcal{C}^{\perp}\,=\, \lbrace v \in \mathbb{F}_{q}^{n} \, | \, (u\cdot v) = 0 \; \; \forall \; u \in \cC \rbrace .
\end{displaymath}
Clearly, $\cC^{\perp}$ is a linear code of dimension $n-\dim(\cC)$. A standard fact from coding theory
asserts that $\cC$ is MDS if and only if $\cC^\perp$ is MDS. 

Now, fix a vector $y \in (\F_q^*)^n$. It defines a {\it $y$-dot product} on $\F_q^n$ given by $(u \cdot v)_{y} = \sum_{i=1}^{n}y_{i} \, u_{i} \, v_{i}$.
If $ y=(1,\dots,1)$, it is the standard dot product. Define 
\begin{displaymath}
\cC^{{\bot}_{y}}= \lbrace v \in \mathbb{F}_{q}^{n} \ | \ (u \cdot v)_{y}=0 \; \; \forall \; u \in \cC \rbrace .
\end{displaymath}
It is easy to see that $\cC^{{\bot}_{y}}$ is equivalent to $\cC^{\perp}$. In fact, $y\,\cC^{\perp_y}=\cC^{\perp}$ in the above notation.

\begin{Def}\label{def:quasi-sd}
A code $\mathcal{C}$ is called \textit{quasi-self-dual} with respect to $ y \in (\mathbb{F}_{q}^{*})^{n} $ if $\mathcal{C} \, = \, \mathcal{C}^{{\bot}_{y}}$.\\
If $ y=(1,1,\dotsc,1)$, we say $ \mathcal{C} $ is \textit{self-dual code}.
\end{Def}

Clearly, if $\cC$ is a quasi-self-dual code with respect to $y=(x_1^2, \dotsc,x_n^2)$, for some $x_i\in\F_q^*$, then 
$x\,\cC$ is a self-dual code. This implies the following proposition. 

\begin{Prop}\label{any quasi-self dual eq sd}
If $char(\F_q)=2$, then any quasi-self-dual code is equivalent to a self-dual code.
\end{Prop}

\subsection{Toric complete intersection codes}\label{S:TCIC}
Recall the definition of a toric complete intersection code following [\ref{ref:Soprunov toric}]. 
Let $\K$ be a field, $\bar{\K}$ be its algebraic closure, and $\K^*=\K\setminus\{0\}$. We use standard terminology and notation
from the theory of Newton polytopes. An element $f$ of the Laurent polynomial ring $ \mathbb{K} \, [t_{1}^{\pm{1}},\dotsc,t_{d}^{\pm{1}}] $
is a finite sum 
\begin{displaymath}
f=\sum_{a \in \cA} c_{a}t^{a}, \text{ where } t^{a}=t_{1}^{a_{1}} \dotsm t_{d}^{a_{d}},\ \cA\subseteq \Z^d,\ c_{a} \in \mathbb{K}.
\end{displaymath} 
The convex hull of the finite set $\cA\subseteq \Z^d$ is called  the {\it Newton polytope} of $f$ and will be denoted by $P(f)$.

For any set $A\subseteq\mathbb{R}^{d}$ let $A_{\mathbb{Z}}\,=\,A\cap\mathbb{Z}^{d}$ denote the set of lattice points in $A$. 
By a slight abuse of notation we use either $|A_\Z|$ or $|A|_\Z$ to denote the cardinality of the set $A_\Z$. 

All polytopes considered in this paper are assumed to be {\it lattice polytopes}, i.e. convex hulls of finitely many points in $\Z^d$.
A polytope of dimension $d$ will be called a $d$-polytope, for short. 

The point-wise sum of two polytopes $P+Q=\{p+q\in\R^d\ |\ p\in P, q\in Q\}$ is called the Minkowski sum. Recall
that any polytope (in fact, any convex body) $P$ is 
uniquely determined by its {\it support function} $l_P$ defined by
$$l_P(v)=\max\{(u\cdot v) \ |\ u \in P\}\quad\text{ for all }v \in \mathbb{R}^d.$$
We will need the following basic properties of the support function:
(1) $l_{P+Q}=l_{P}+ l_{Q}$ and (2) $P \subseteq Q$ if and only if $l_P(v) \leq l_Q(v)$ for every $v$ in $\mathbb{R}^d$.

We denote the Euclidean $d$-dimensional volume of $P$
by $V_d(P)$, or simply by $V(P)$ when the dimension is clear. We use $V(P_1,\dots, P_d)$ to denote the normalized 
{\it mixed volume} of $d$ lattice polytopes $P_1,\dots, P_d$. By definition, 
$$V(P_1,\dots, P_d)=\sum_{I\subseteq \{1,\dots,d\}} (-1)^{d-|I|}V_d(P_I),$$
where $P_I=\sum_{i\in I}P_i$. The mixed volume $V(P_1,\dots, P_d)$ is non-negative, 
multilinear with respect to Minkowski addition, and coincides with $d!V_d(P)$ when $P_i=P$ for every $1\leq i\leq d$. 
More about  the mixed volume can be found in \cite[Ch. 4]{Geomineq}.

Now fix a finite subset $S\,=\,\lbrace p_{1},\dotsc,p_{n} \rbrace $ of the algebraic torus ${(\mathbb{K}^*)^{d}} $ and a finite-dimensional subspace 
$\mathcal{L}$ of $ \mathbb{K} \, [t_{1}^{\pm{1}},\dotsc,t_{d}^{\pm{1}}] $.

\begin{Def}\label{D:evalcode}
Define the \textit{evaluation map}
\begin{displaymath}
ev_{S} \, : \, \mathcal{L}\rightarrow \mathbb{K}^n, \; f \, \mapsto \, (f(p_{1}),\dotsc,f(p_{n})).
\end{displaymath}
The image of $ev_{S}$ is called the \textit{evaluation code} corresponding to $S$ and $\cL$. 
We will denote this code by $ \mathcal{C}_{S,\mathcal{L}}$.
\end{Def}

Clearly, $ \mathcal{C}_{S,\mathcal{L}}$ is a linear code over $\K$ of block-length $n$.

Toric complete intersection codes are special evaluation codes when $S$ is the solution set of a Laurent polynomial system satisfying some assumptions. Here is the precise definition.
\begin{Def}\label{D:TCI}
Fix a collection of $d$-polytopes $P_{1},\dotsc,P_{d}$ in $\mathbb{R}^{d}$ and consider $d$ Laurent polynomials $f_{1},\dots,f_{d}$ over $\mathbb{K}$ with Newton polytopes $P_{1},\dots,P_{d}$ such that the solution set $S$ of the system $f_{1}=\dotsb=f_{d}=0$ in $({\bar{\K}}^*)^d$ satisfies the following:
\begin{enumerate}
\item[\textit{(1)}] $|S|=V(P_1,\dotsc ,P_d)$,
\item[\textit{(2)}] the set $S$ consists of $\mathbb{K}$-rational points i.e. $S \subseteq ({{\K}}^*)^d$.
\end{enumerate}
Then $S$ is called a {\it toric complete intersection over $\K$}.
\end{Def}

\begin{Rem}
In general, the set $S$ is the intersection of $n$ hypersurfaces in a toric variety associated with the polytope $P$.
According to the Bernstein-Kushnirenko-Khovanskii bound \cite{B, Kush}, 
if $S$ consists of isolated points, its cardinality $|S|$ cannot
exceed the mixed volume $V(P_1,\dotsc ,P_d)$. Moreover, the bound is attained for systems with generic coefficients
(having the $P_i$ fixed) in which case the hypersurfaces do not intersect outside of the torus $(\bar\K^*)^d$ and the intersections are transversal.
This is guaranteed by the assumption {\it (1)}. In particular, this implies that the local intersection multiplicities equal one, and the ideal $\langle f_1,\dots,f_d\rangle$ is radical. 
\end{Rem}

The following toric analog of the Euler--Jacobi theorem by Khovanskii \cite{uspehi} 
is fundamental for our results about toric complete intersection codes and their duals.
For a proof that works over arbitrary algebraically closed fields see \cite[Sec. 14]{Kunz}.
First we need the following definition.

\begin{Def}
Let $f_{1},\dotsc,f_{d} \in \mathbb{K} \, [t_{1}^{\pm{1}},\dotsc,t_{d}^{\pm{1}}]$ be Laurent polynomials. The Laurent polynomial 
\begin{displaymath}
J^{\mathbb{T}}_{f}\,=\det\bigg(t_{j}\dfrac{\partial f_{i}}{\partial t_{j}}\bigg)
\end{displaymath}
is called the \textit{toric Jacobian} of $f_{1},\dotsc,f_{d}$.
\end{Def}

\begin{Th}\cite{uspehi} \label{T:thJacobian}
Let $S$ be a toric complete intersection over $\K$. 
Let $ P = P_{1}+\dotsb+P_{d}$ be the Minkowski sum and $P^\circ$ be its interior. Then for any $h \in \mathcal{L}(P^{\circ})$ we have 
$$\sum_{p \in S}\frac{h(p)}{J^{\mathbb{T}}_{f}(p)}=0.$$
\end{Th} 

Note, since the local intersection multiplicities equal to one, $J^{\mathbb{T}}_{f}(p) \neq 0 $ for every $ p \in S$, and 
the above sum makes sense.

To finish the definition of the toric complete intersection code we describe the space $\cL$. 
As before, let $P^\circ$ be the interior of  $P=P_1+\dots+P_d$.
Fix any subset $A$ of $P^{\circ}$. It defines a space of Laurent polynomials over $\K$:
\begin{displaymath}
\mathcal{L}(A)\,=\spn_{\mathbb{K}}\lbrace t^{a}\,|\,a \in A_{\mathbb{Z}} \rbrace \, \subseteq \mathbb{K} \, [t_{1}^{\pm{1}},\dotsc,t_{d}^{\pm{1}}].
\end{displaymath}

\begin{Rem}\label{R:incube}
 It is clear that $\dim_{\bar\K}\cL(A)$ equals the cardinality of $A_\Z$, but $\dim_{\K}\cL(A)$ may be smaller.
When $\K=\F_q$ a finite  field, $\dim_{\bar\K}\cL(A)=\dim_{\K}\cL(A)$ if and only if different points of $A_\Z$ represent 
different classes in $\Z^d/((q-1)\Z)^d$. In particular, this happens if $P$ is contained in the $d$-cube $([0,q-2])^d$.
Loosely speaking, this means that the degrees of the monomials appearing in the system $f_1=\dots=f_d$ are 
small compared to the size of the field. In what follows we will always make this assumption and write $\dim$ for
$\dim_\K$ or $\dim_{\bar\K}$.
\end{Rem}

\begin{Def}
Let $S$ be toric complete intersection over $\K$. Let $A\subseteq P^{\circ}$ and let $\cL(A)$ be the corresponding
polynomial space. The evaluation code $\cC_{S,{\cL}(A)}$ is called a \textit{toric complete intersection code}, denoted simply by 
$\mathcal{C}_{S,A}$.    
\end{Def}

In \cite{So1} the second author gave lower bounds for the minimum distance of toric complete intersection codes. It turns out
that the bound is significantly better if the solution set $S$ satisfies an extra assumption of ``generic position". We formulate it below.

\begin{Def}\label{D:(3)} A subset $S\subset (\K^*)^d$ is said to be in {\it $Q$-generic position} if
there exists a $d$-polytope $Q$ such that for any subset $T \subseteq S$ of size $|Q_\Z|$ the evaluation map 
$ev_T: \cL(Q) \rightarrow \K^{|Q_{\Z}|}$ is an isomorphism.
\end{Def}

In other words, $S$ is in $Q$-generic position if for any collection $T$ of size $|Q_\Z|$ there is a polynomial
$h\in \cL(Q)$ which takes the zero value at all but the last point of $T$. 
For example, when $Q=\triangle_d$ is the standard $d$-simplex, i.e. 
the convex hull of $\{0,e_1,\dots, e_d\}$, where $\{e_1,\dots, e_d\}$ is the standard basis for $\R^d$,
this means that no $d+1$ points of $S$ lie on a hyperplane.

Here is the lower bound on the minimum distance for toric complete intersection codes.
\begin{Th}\cite{So1}\label{T:bound}
Let $S$ be a toric complete intersection in $Q$-generic position.
Let $A$ be any set such that $A+mQ \subseteq P^{\circ}$ up to a lattice translation, for some $m \geq 0$. Then
$$
d(\cC_{S,A}) \, \geq \, (|Q_{\Z}|-1)m+2.
$$
\end{Th}

We end this section with geometric conditions on the polytopes $P_1,\dots,P_d$ which insure that assumption (3) is satisfied for
generic systems. Intuitively, it says that $Q$ is small compared to the $P_i$.

\begin{Th}\cite{So1}\label{T:cond-on-poly}
Let $Q$ be a $d$-polytope such that $Q_{\Z}$ generates $\Z^d$. Suppose
\begin{enumerate}
\item[a.] $V(P_1,\dotsc ,P_{d-1},Q) \, \geq \, |Q_{\Z}|$,
\item[b.] $(|Q_{\Z}|-1)Q \subseteq P_d$, up to a lattice translation, or $P_1,\dots, P_d$ and $Q$ have the same normal fan and $P_1+\dots+P_{d-1}+Q\subseteq P_d$, up to a lattice translation. 
\end{enumerate}
Then the solution set of any system $f_{1}=\dotsb=f_{d}=0$ with Newton polytopes $P_1,\dotsc ,P_d$ and generic coefficients 
is in $Q$-generic position.
\end{Th}

\section{Results in arbitrary dimension $d$} 



We begin this section with an immediate result from the assumption \textit{(3)}. 

\begin{Th}\label{T:CQ is MDS} Let $S\subseteq (\K^*)^d$ be any subset in $Q$-generic position for some $d$-polytope $Q$.
Then the evaluation code $\cC_{S,Q}$ is an MDS code.
\end{Th}

\begin{pf} Denote $\cC:= \cC_{S,Q}$.
We need to show that $\cC$ is an $[n,k,n-k+1]_{q}$-code where $k=\dim(\cC)$ and $n=|S|$. 
First, we show that $k=|Q_{\mathbb{Z}}|$. Consider the evaluation map
\begin{displaymath}
ev_{S}\,:\,\mathcal{L}(Q) \rightarrow \mathbb{K}^{n},\quad f\mapsto(f(p_{1}),\dotsc , f(p_{n})).
\end{displaymath} 
By definition $\cC=\Img(ev_S)$. Since $\dim\cL(Q)=|Q_\Z|$, it is enough to show that $ev_S$ is injective.
If $f\in\Ker(ev_S)$ then $f\in\Ker(ev_T)$ for any subset $T\subseteq S$ of size $|Q_{\mathbb{Z}}|$.
By \rd{(3)}, $ev_{T}$ is an isomorphism, so $\Ker(ev_T)$ is trivial. Therefore $f=0$. 

Now we show that $d(\cC)=n-k+1$. By before, $\Ker(ev_T)$ is trivial for any
$T\subseteq S$ of size $|Q_{\mathbb{Z}}|$. Therefore any non-zero $f\in\cL(Q)$ can have at most 
$|Q_\Z|-1$ zeroes in $S$. In other words, the image of $f$ under  $ev_S$ has weight
at least $n-|Q_\Z|+1$. This shows that $d(\cC)\geq n-k+1$. On the other hand,
by the Singleton bound $d(\cC)\leq n-k+1$. This proves that $\cC$ is an $[n,k,n-k+1]_{q}$-code.
\end{pf}

\begin{Cor}
The dual code $\cC_{S,Q}^{\perp_{y}}$ is an MDS code.
\end{Cor}

This follows from the fact that $\cC^{\perp_{y}}$ is equivalent to $\cC^\perp$ and the dual of an MDS-code is
also MDS, as we mentioned in the Preliminaries.

The following theorem relates the toric complete intersection codes defined by $A\subseteq P^\circ$ and $B\subseteq P^\circ$
which satisfy $A+B\subseteq P^\circ$. Here $A+B=\{a+b\in\R^d\ |\ a\in A, b\in B\}$ is the Minkowski sum.

\begin{Th}\label{T:y-dual toric} Let $S$ be a toric complete intersection.
Let  $A, B$ be subsets of $P^\circ$ such that $A+B\subseteq P^\circ$. If
$\dim(\cC_{S,B})=|S|-\dim(\mathcal{C}_{S,A})$, then there exists $y\in(\K^*)^n$ such that
$$\cC_{S,B}=\cC_{S,A}^{\perp_{y}}.$$
In particular, if  $|S|$ is even, $2A \subseteq P^{\circ}$, and $\dim(\cC_{S,A})=|S|/2$ then $\cC_{S,A}$ is  {quasi-self-dual}. 
\end{Th}

\begin{pf}  Let $S=\{p_1,\dots, p_n\}$.
First, for any $h=fg$, where $f \in \mathcal{L}(A)$ and $ g \in \mathcal{L}(B)$, we have $h \in \mathcal{L}(A+B)\subseteq \mathcal{L}(P^{\circ})$. 
By \rt{thJacobian}, 
\begin{displaymath}
\sum_{i=1}^n\frac{h(p_i)}{J^{\mathbb{T}}_{f}(p_i)} = \sum_{i=1}^n\frac{f(p_i)g(p_i)}{J^{\mathbb{T}}_{f}(p_i)} = 0. 
\end{displaymath}
This implies that  $ev_S(f)$ and $ev_S(g)$ are $y$-orthogonal where 
$y\,=\,\Big(\dfrac{1}{J^{\mathbb{T}}_{f}(p_{1})},\dotsc,\dfrac{1}{J^{\mathbb{T}}_{f}(p_{n})}\Big)$. 
Hence $\cC_{S,B}$ is a subspace of $\cC_{S,A}^{\perp_{y}}$. On the other hand, 
$\dim(\cC_{S,B})=n-\dim(\cC_{S,A})=\dim\big(\cC_{S,A}^{\perp_{y}}\big)$, and the first statement follows. 

For the second part, let $B=A$. Then, $\dim(\cC_{S,A})=\dim(\cC_{S,B}) =n/2$. By Definition \ref{def:quasi-sd}, $\cC_{S,A}$ is a quasi-self-dual code.    
\end{pf}

As an immediate consequence of the above theorem and Proposition \ref{any quasi-self dual eq sd} we obtain the following.

\begin{Cor}\label{C:CA eq self dual}
Let $char(\K)=2$. If  $|S|$ is even, $2A \subseteq P^{\circ}$, and $\dim(\cC_{S,A})=|S|/2$ then $\cC_{S,A}$ is equivalent to a self-dual code.
\end{Cor}

Our ultimate goal is to give a description of the polytopes $P_1,\dots, P_d$ and the set $A$ such that generic systems 
produce quasi-self-dual codes. For this we need a way to compute the dimension $\dim(\cC_{S,A})$. According to \rd{evalcode},
this amounts to computing the dimension of the kernel of the evaluation map, 
\begin{equation}\label{e:dim of C_A}
 \dim(\cC_{S,A})=\dim(\cL(A))-\dim\Ker(ev_S)=|A_\Z|-\dim\Ker(ev_S). 
\end{equation}
Let $J=\langle f_1,\dots,f_d\rangle$ be radical. Then polynomials in $\Ker(ev_S)$ are, in fact, elements of $\cL(A)\cap J$.
In other words, one has to compute an analog of the Hilbert function for the ideal~$J$:
$$\Hilb_J(A)=\dim(\cL(A)\cap J).$$
Although this can be done in some situation, there appears to be no simple formula for $\Hilb_J(A)$ in general.
We explore the $d=2$ case in the next section. Also, in \cite{SaSo} the authors give a formula 
for $\dim(\cC_{S,A})$ when the polynomials $(f_1,\dots, f_d)$ give rise to 
a regular sequence $(F_1,\dots, F_d)$ in the {\it homogeneous coordinate ring} of a toric variety. 
We plan to return to this problem in the future.

\section{Results in dimension $d=2$} 

In this section we concentrate on the case $d=2$. We reserve the word ``polygon" for 
any convex polytope of dimension at most two.
Let $S\subseteq(\K^*)^2$ be a toric complete intersection defined by
Laurent polynomial system $f_1=f_2=0$ with lattice polygons $P_1, P_2$ as in \rd{TCI}.
As before, $P^\circ$ denotes the interior of $P=P_1+P_2$, and $V(P,Q)$
the normalized mixed volume (mixed area) of $P$ and $Q$, i.e.
$$V(P,Q)=V(P+Q)-V(P)-V(Q),$$
where $V(P)$ is the Euclidean area of $P$. It is easy to check that $V(P,Q)=0$
if and only if either one of the polygons is a point or $P$, $Q$ are parallel segments.

Our goal is to give a description of lattice polygons $P_1$, $P_2$ for which there exists $A$
satisfying 
\begin{equation}\label{e:problem}
2A\subseteq (P_1+P_2)^\circ,\quad\text{and}\quad \dim(\cC_{S,A})=V(P_1,P_2)/2.
\end{equation}
Then,  $\cC_{S,A}$ is quasi-self-dual, by \rt{y-dual toric}. 

In \rt{geom criterion} below we give a general geometric condition on $P_1$, $P_2$, and $A$
that guarantees that $\cC_{S,A}$ is quasi-self-dual. Then we  look at special cases (\rp{unmixed}, \rt{summands})
when we can construct $P_1$, $P_2$, and $A$ explicitly.

Intuitively, $A$ has to be just a bit ``smaller'' than the ``average'' of $P_1$ and $P_2$. Although
we have Minkowski addition on the space of lattice polygons, there is no subtraction, in general.
To resolve this, we introduce the following analog of difference of (convex) sets.

Let $A$, $B$ we subsets of $\R^d$. Define
$$A-B=\{u\in \R^d\ |\ u+B\subseteq A\}.$$

It is easy to show that $A-B$ is convex if $A$ is convex. Also $(A-B)+B\subseteq A$, but not equal to $A$, in general.
Rather, it is the largest subset in $A$ that has $B$ as a Minkowski summand.

Now we are ready to give a combinatorial formula for $\dim(\cC_{S,A})$. We begin with a few lemmas. 

\begin{Lemma}\label{L:poly-compare}
Let $A\subset \R^2$. Then $V = \langle f_{1} \rangle \cap \cL(A)$ is a subspace
of $\cL(A)$ with a basis $\cB=\lbrace f_{1}t^{a} \, | \, a \in (A-P_1)_{\mathbb{Z}} \rbrace$.  
\end{Lemma}

\begin{pf} The fact that $V\subseteq\cL(A)$ is a subspace is straightforward. Denote $R=A-P_1$.
To show that $\cB$ is linearly independent, suppose
 $${\sum}_{a \in R_{\Z}}\lambda_{a}f_{1}t^{a}=0,\quad \lambda_{a} \in \mathbb{K}.$$
 Then
$ f_{1}\left(\sum_{a \in R_{\Z}}\lambda_{a}t^{a}\right)=0$ in $\cL(A)$.
Since $\mathcal{L}(A)$ is a subset of $\mathbb{K}[t_1^{\pm 1},t_2^{\pm 1}]$, it has no zero divisors. 
Thus,  $\sum\lambda_{a}t^{a}=0$, which implies that $\lambda_{a}=0$, and so $\cB$ is linearly independent.

To show $\cB$ spans $V$, note that any $g \in V$ can be written as $g=hf_{1} \in \mathcal{L}(A)$. We have
$P(g)\subseteq A$ and $P(g)=P(h)+P_{1}$. Thus $P(h)\subseteq A-P_1=R$. In particular, every monomial in $h$
has exponent lying in $R_\Z$, i.e. $h= \sum_{a \in R_{\Z}} \lambda_{a}t^{a}$. This shows that $g$ is a linear
combination of elements in $\cB$.
\end{pf}

\begin{Lemma} \label{equality for kernel} Let $f_1$ be absolutely irreducible and  $A$ a lattice polygon. If $V(P_1,A) < V(P_1,P_2)$ then
\begin{displaymath}
 \la f_1,f_2 \ra \cap \: \cL(A)=\la f_1\ra \cap\: \cL(A). 
 \end{displaymath}
\end{Lemma}

\begin{pf}
One inclusion $\la f_1 \ra \cap \: \cL(A) \subseteq \la f_1,f_2  \ra \cap \: \cL(A)$ is obvious. For the other one, consider $f \in \la f_1,f_2\ra \cap \: \cL(A)$. Clearly, $f$ vanishes at points in $S$. Now,  the system $f_1=f=0$ has at least $|S|=V(P_1,P_2) > V(P_1,A)$ solutions. On the other hand, $f \in \cL(A)$ implies $P(f) \subseteq A$, hence, by the Bernstein-Kushnirenko theorem, $f$ and $f_1$ must have a common component. Since $f_1$ is absolutely irreducible, $f_1$ divides $f$.  Therefore, $f \in \la f_1 \ra \cap \: \cL(A)$. This implies $\la f_1,f_2\ra\cap \: \cL(A) \subseteq \la f_1\ra \cap \cL(A)$, and the statement follows.
\end{pf}

\begin{Prop} \label{P:dim}
Let $S$ be a toric complete intersection over $\K$
and suppose $f_1$ is absolutely irreducible.
Let $A$ be a lattice polygon  such that $V(P_1,A)< V(P_1,P_2)$. Then
  $$\dim(\cC_{S,A})=|A|_\Z-|A-P_1|_\Z.$$
\end{Prop}

\begin{pf}
By \re{dim of C_A} we have $\dim(\cC_{S,A})=|A_\Z|-\dim\Ker(ev_S)$.
Since $\la f_1, f_2\ra$ is radical, it implies that  $\Ker(ev_S)=\la f_1, f_2\ra \cap \cL(A)$. The latter equals
$\la f_1\ra \cap \cL(A)$, by Lemma \ref{equality for kernel}.  The result now follows from  \rl{poly-compare}.
\end{pf}

\begin{Rem}\label{R:dim}
Note that $A=P_1$ corresponds to $A-P_1={(0,0)}$, the origin. In this case, $\dim(\cC_{S,A})=|A|_\Z-1$.
If $A$ does not contain any lattice translate of $P_1$ then $A-P_1$ is empty and $\dim(\cC_{S,A})=|A|_\Z$.
\end{Rem}

Now \rp{dim} and \rt{y-dual toric} provide the following geometric criterion.

\begin{Th}\label{T:geom criterion}
Let $S$ be a toric complete intersection over $\K$
and suppose $f_1$ is absolutely irreducible. Let $A$ be a lattice polygon  such that
\begin{enumerate}
\item[i.] $V(P_1,A)< V(P_1,P_2)$,
\item[ii.]  $2A\subseteq (P_1+P_2)^\circ$,
\item[iii.] $|A|_\Z-|A-P_1|_\Z=V(P_1,P_2)/2$.
\end{enumerate}
Then $\cC_{S,A}$ is a quasi-self-dual toric complete intersection code.
\end{Th}

To make our result more explicit we analyze the geometric conditions of \rt{geom criterion}
in some special cases. First we
consider the so-called {\it unmixed} case, when $P_1$ and  $P_2$ are integer dilates
 of the same lattice polygon $Q$. In other words, $P_1=m_1Q$ and $P_2=m_2Q$, for some positive integers $m_i$.
 Choose $A=aQ$ for some positive integer $a$. Then the conditions in  \rt{geom criterion}
 become
 \begin{equation}\label{e:geomproblem unmixed}
a<m_2,\quad 2a < m_1+m_2,\quad\text{and}\quad |aQ|_\Z-|(a-m_1)Q|_\Z=m_1m_2V(Q).
\end{equation}
We have the following result.

\begin{Prop}\label{P:unmixed} 
Let $P_1=m_1Q$, $P_2=m_2Q$, and $A=aQ$ for some lattice polygon $Q$ and positive integers $m_1$, $m_2$, and $a$.
Suppose \re{geomproblem unmixed} holds. Then only the following three cases are possible.
\begin{enumerate}
\item $Q$ is $\GL(2,\Z)$-equivalent to the standard 2-simplex, $a=(m_1+m_2 - 3)/2$, and $a\in\N$;
\item $Q$ is $\GL(2,\Z)$-equivalent to either the triangle with vertices $\{0, 2e_1, e_2\}$ or the standard square,  
$a=(m_1+m_2 - 2)/2$, and $a\in\N$;
\item $Q$ is $\GL(2,\Z)$-equivalent to one of the sixteen Fano polygons in \rf{fano}, $a=(m_1+m_2 - 1)/2$, and $a\in\N$.
\end{enumerate}
\end{Prop}

 \begin{figure}[h]
\includegraphics[scale=.7]{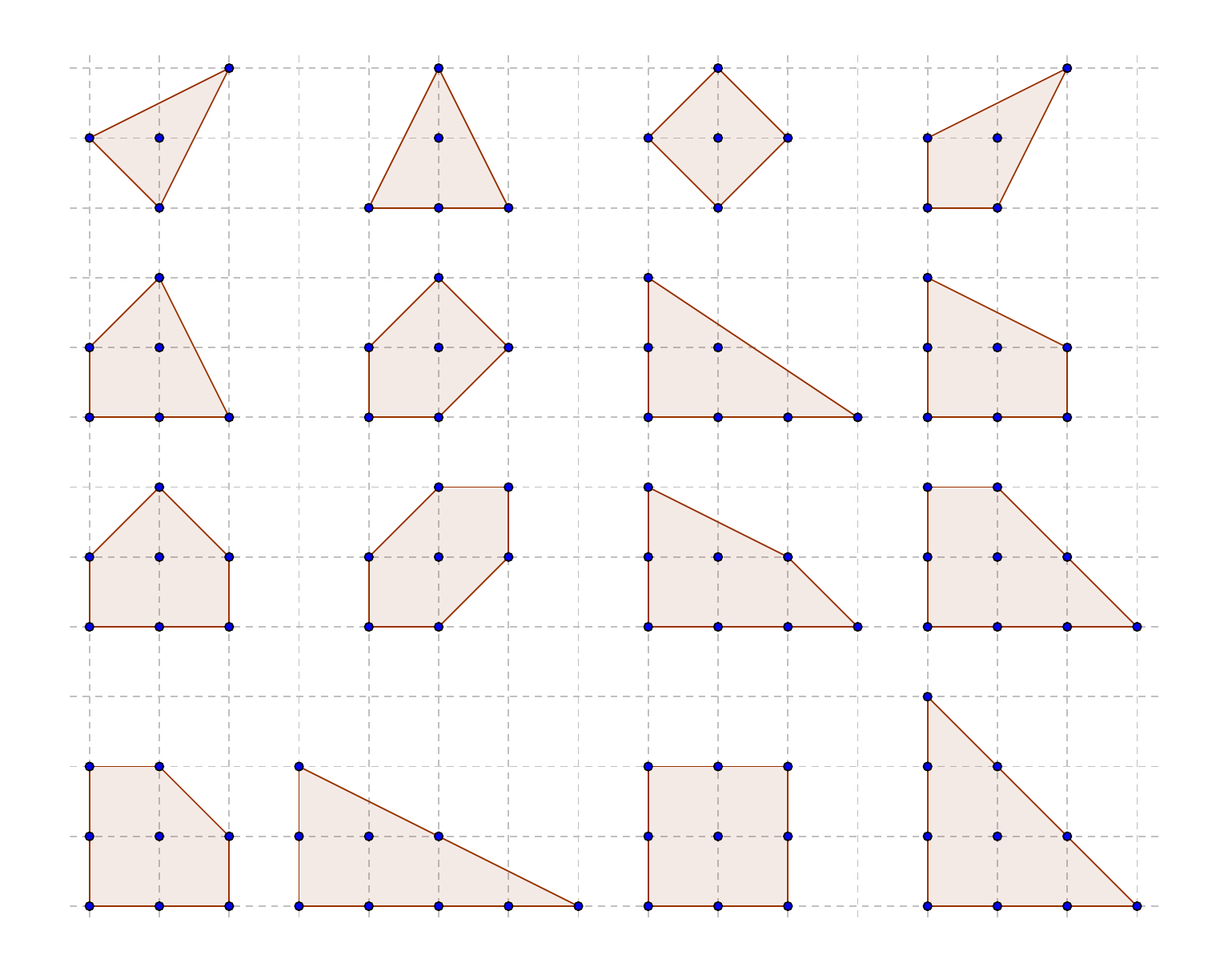}
\caption{The sixteen $\GL(2,\Z)$-classes of Fano polygons.}
\label{F:fano}
\end{figure}

\begin{pf} According to Pick's formula $|aQ|_\Z=a^2V(Q)+\frac{a}{2}|\partial Q|_\Z+1$, where $\partial Q$ denotes the boundary of $Q$.
This is the Ehrhart polynomial of $Q$.  
 \rf{Pick} depicts the set of all $(c_1,c_2)$ (marked with dots) which are possible coefficients of Ehrhart polynomials,
 i.e. for which there exists a lattice polygon $Q$ with $c_1=\frac{1}{2}|\partial Q|_\Z$ and $c_2=V(Q)$, see \cite{Beck}. 
 
  \begin{figure}[h]
\includegraphics[scale=.7]{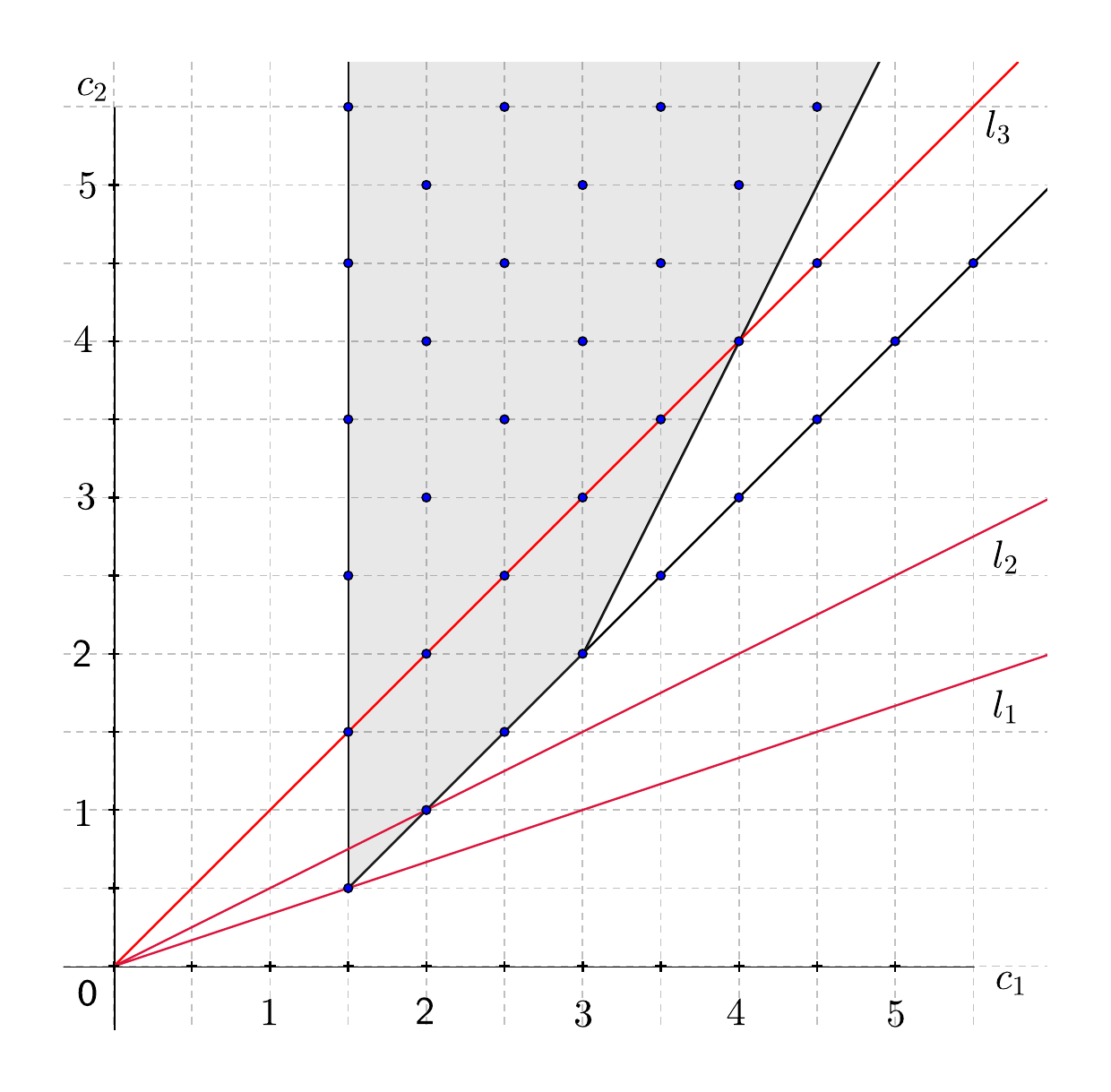}
\caption{The set of coefficients of Ehrhart polynomials.}
\label{F:Pick}
\end{figure}

 These points $(c_1,c_2)$ have integer or  half-integer  coordinates and 
 consist of points lying either in the shaded region or on the line $c_2=c_1-1$ with the exception of a single point $(9/2,9/2)$.

First, assume $a\geq m_1$. Applying Pick's formula to $|aQ|_\Z$ and $|(a-m_1)Q|_\Z$ and simplifying, we 
see that the equation in \re{geomproblem unmixed} is equivalent to
$$(m_1+m_2-2a)V(Q)=\frac{1}{2}|\partial Q|_\Z.$$
 It follows from  \rf{Pick} that the only lines $\lambda c_2= c_1$ with $\lambda\in\N$ that intersect the set of possible coefficients 
 are $3c_2= c_1$, $2c_2= c_1$, and  $c_2= c_1$, labelled by $l_1$, $l_2$, and $l_3$, respectively.
 
 In the first case, $c_1=3/2$, $c_2=1/2$, which corresponds to $Q$ being $\GL(2,\Z)$-equivalent to the standard 2-simplex.
In this case $a=(m_1+m_2 - 3)/2$ and it has to be a positive integer. 
In the second case, $c_1=2$, $c_2=1$, which corresponds to $Q$ being $\GL(2,\Z)$-equivalent to
either the triangle with vertices $\{0, 2e_1, e_2\}$ or the standard square. Here $a=(m_1+m_2 - 2)/2$, and we must have $a\in\N$.
Finally, $c_2=c_1$ corresponds to lattice polygons with exactly one interior lattice point. These are Fano polygons and
there are exactly sixteen classes of them up to $\GL(2,\Z)$ equivalence. In this case $a=(m_1+m_2 - 1)/2$, and $a$ must be in $\N$.

Now assume  $a< m_1$.  In this case $\dim(\cC_{S,A})=|A_\Z|$ by \rr{dim}, and the equation in \re{geomproblem unmixed}
becomes $|aQ|_\Z=m_1m_2V(Q)$. Again, by using Pick's formula one can show that this is equivalent to
the line $(m_1m_2-a^2)c_2=ac_1+1$ having a non-trivial intersection with the set of lattice points in \rf{Pick}, which is impossible if $1\leq a<m_1$. 
This completes the proof of the theorem.
\end{pf}

\begin{Rem} We point out that the last three polygons in the bottom row in \rf{fano} are, in fact, particular cases
of (2) when both $m_1$ and $m_2$ even, and (3) when both $m_1$ and $m_2$ are multiples of three.
\end{Rem}

Combining the results of \rt{bound},  \rt{geom criterion}, and \rp{unmixed} we obtain the following.

\begin{Cor}\label{C:unmixed} 
Let $S$ be a toric complete intersection over $\K$ in $Q$-generic position and assume $f_1$ is absolutely irreducible.
Let $P_1=m_1Q$ , $P_2=m_2Q$ and $A=aQ$ be as (1), (2) or (3) in \rp{unmixed}.
Then $\cC_{S,A}$ is quasi-self-dual with parameters 
\begin{enumerate}
\item $n=m_1m_2$, $k=n/2$, $d(\cC_{S,A}) \geq(m_1+m_2+1)/2$; or
\item $n=2m_1m_2$, $k=n/2$, $d(\cC_{S,A}) \geq m_1+m_2$; or
\item $n=2V(Q)m_1m_2$, $k=n/2$, $d(\cC_{S,A})\geq V(Q)(m_1+m_2-1)+2$,
\end{enumerate}
respectively.
\end{Cor}

Notice that the geometric conditions a.--b. in \rt{cond-on-poly} hold for $m_1Q$ , $m_2Q$,
as long as $1<m_1<m_2$. Therefore, systems $f_1=f_2=0$ with Newton
polygons $m_1Q$ , $m_2Q$ and {\it generic coefficients}  in $\bar\K$ will
produce quasi-self dual codes $\cC_{S,A}$.

\medskip

Our next situation is more general. Here we only assume that $P_1$ is a Minkowski
summand of $A$,  and $A$ is a Minkowski summand of $P_2$. In other words,
$$A=P_1+R_1\quad\text{and}\quad P_2=A+R_2,$$  
for some lattice polygons $R_1$, $R_2$ (we allow $R_1$ to be a point or a lattice segment).

Recall that $l_P(v)$ denotes the support function of $P$.
Also, by $\Fan(P)$ we mean the set of primitive 
lattice vectors (i.e. whose entries are coprime)  that are the outer normals to the edges of $P$.

\begin{Th}\label{T:summands} 
Let $A=P_1+R_1$ and $P_2=A+R_2$ for some lattice polygons $P_1$, $R_1$, and $R_2$. 
Then i.--iii. in \rt{geom criterion} hold if and only if $R_1\subset R_2^\circ$ and
$l_{R_2}(v)= l_{R_1}(v)+1$ for all  $v\in\Fan(P_1)$.
\end{Th}

\begin{pf} The condition $2A\subseteq (P_1+P_2)^\circ$ written in terms of the support functions
translates to $l_{2A}(v) < l_{P_1+P_2}(v)$ for all $v\in \R^2$. By properties of the support function
this is equivalent to $l_{R_1}(v) < l_{R_2}(v)$ for all $v\in \R^2$, which means $R_1\subset R_2^\circ$.

Next  we look at condition iii:
\begin{equation}\label{e:iii}
2|A|_\Z-2|A-P_1|_\Z=V(P_1,P_2).
\end{equation}
 
Applying Pick's formula and linearity of the mixed volume, we can rewrite the left hand side as follows.
$$2|P_1+R_1|_\Z-2|R_1|_\Z=2V(P_1+R_1)-2V(R_1)+|\partial P_1|_\Z=V(P_1,P_1+2R_1)+|\partial P_1|_\Z,$$
where we used an obvious relation $|\partial(P_1+R_1)|_\Z=|\partial P_1|_\Z+|\partial R_1|_\Z$.
Now \re{iii} is equivalent to
\begin{equation}\label{e:iii equiv}
V(P_1,R_1)+|\partial P_1|_\Z=V(P_1,R_2).
\end{equation}
There is an ``inductive'' formula  for computing the mixed volume \cite[Ch. 4]{Geomineq}.
It can be adapted to the lattice situation. In dimension two it states the following. 
Let $P$ be a lattice polygon and $L_v$ be the lattice length of 
the edge of $P$ corresponding to $v\in \Fan(P)$. Then for any lattice polygon $R$ 
$$V(P,R)=\sum_{v\in\Fan(P)}l_R(v)L_v$$
Note that when all $l_R(v)$ equal one, the above sum 
is just $|\partial P_1|_\Z$. Therefore, \re{iii equiv} is equivalent to
\begin{equation}\label{e:iii equiv equiv}
\sum_{v\in\Fan(P_1)}(l_{R_1}(v)+1)L_v=\sum_{v\in\Fan(P_1)}l_{R_2}(v)L_v.
\end{equation}
On the other hand, we have $l_{R_1}(v) < l_{R_2}(v)$ for all $v\in \R^2$, by condition ii.
In particular, $l_{R_1}(v)+1 \leq l_{R_2}(v)$ for $v\in\Fan(P_1)$, as $l_{R_i}(v)$ takes 
integer values for these $v$. Therefore \re{iii equiv equiv} holds if and only if
$$l_{R_1}(v)+1= l_{R_2}(v)\quad\text{for all } v\in\Fan(P_1).$$

Finally, the condition $V(P_1,A)<V(P_1,P_2)$ is the same as $V(P_1,R_2)>0$, which is true
since $P_1$ is 2-dimensional and $R_2$ is not just a point, otherwise $R_1\subset R_2^\circ$ would be false.

\end{pf}

\begin{Rem}\label{R:general} In fact, we can restate the condition in \rt{summands}
as follows: 
$$2A\subset (P_1+P_2)^\circ\quad\text{and}\quad 2l_A(v)=l_{P_1}(v)+l_{P_2}(v)-1\text{ for all }v\in\Fan(P_1).$$
For this it is enough to only assume that $P_1$ is a Minkowski summand of $A$.
This justifies what we said previously that $A$ has to be a bit smaller than the average of $P_1$ and $P_2$.
However, this condition is not as convenient for constructing examples as the one in  \rt{summands}. 
\end{Rem}

\section{Algorithm and Examples} \label{algorithm}

In this section we collect examples of toric complete intersection codes. All our examples were
produced using MAGMA algebra system \cite{Magma}. Our method is a rather straightforward 
random search for toric complete intersections. The algorithm which we put below works well
for small polygons $P_1$. In all our examples we work over $\F_{16}$. This is the smallest
field of characteristic two for which all the considered polygons lie in the square $[0,q-2]^2$
(see \rr{incube}).

First we need a simple necessary
condition for $S$ to be a toric complete intersection.

\begin{Prop} \label{P:cond-rank}
Let $S$ be a toric complete intersection with Newton polygons $P_1$, $P_2$. 
Then the rank of the evaluation map $ev_S$ satisfies
$$rk(ev_S) \leq |P_2|_\Z-|P_2-P_1|_{\Z}-1.$$
In particular, when $P_1$ is a Minkowski summand of $P_2$ we have
$$rk(ev_S)\leq |S|-|P_1^{\circ}|_\Z.$$
\end{Prop}

\begin{pf} Let $n=|S|=V(P_1,P_2)$ and consider the following sequence which is exact in the first two terms:
\begin{displaymath}
0 \, \rightarrow \, \Ker(ev_{S}) \, \rightarrow \, \mathcal{L}(P_{2}) \, \xrightarrow{\text {$ev_S$}} \, \mathbb{K}^{n}.
\end{displaymath}
Clearly, $\langle f_1\rangle\cap\mathcal{L}(P_2)$ is a subspace of $\Ker(ev_S)$.
On the other hand $f_2$ lies in $\Ker(ev_S)$ and has no common factors with $f_1$, so the inclusion is strict.
Therefore, by \rl{poly-compare}, we have
$$|P_2-P_1|_\Z=\dim \langle f_1\rangle\cap\mathcal{L}(P_2)<\dim \Ker(ev_S)=|P_2|_\Z-rk(ev_S),$$
and the first inequality follows.

Now if $P_2=P_1+R$ for some lattice polygon $R$ then applying Pick's formula,
$$|P_2|_\Z-|P_2-P_1|_{\Z}-1=|P_1+R|_\Z-|R|_\Z-1=V(P_1,R)+V(P_1)+\frac{1}{2}|\partial P_1|_Z-1.$$
By the linearity of the mixed volume $V(P_1,P_2)=V(P_1,P_1+R)=2V(P_1)+V(P_1,R)$, so we get
$$|P_2|_\Z-|P_2-P_1|_{\Z}-1=V(P_1,P_2)-\Big(V(P_1)-\frac{1}{2}|\partial P_1|_Z+1\Big).$$
By Pick's formula again, the expression in the parentheses on the right is $|P_1^\circ|_\Z$.
\end{pf}

Below is the algorithm we use to produce examples of toric complete intersections $S$.
The input is lattice polygons $P_1$, $P_2$, and $Q$ if we wish $S$ to be in $Q$-generic position. 
The output is $S$ and the polynomials $f_1$, $f_2$.

\medskip

\noindent {\bf Algorithm.}

\begin{itemize}
\item[1.] Choose a random absolutely irreducible Laurent polynomial $f_1$ whose Newton polytope is $P_1$.
\item[2.] Find the $\K$-rational points of $f_1=0$.
\item[3.] Choose a subset $S$ of $n=V(P_1,P_2)$ of the points in Step 2  in $Q$-generic position. 
\item[4.] Check whether the rank of the evaluation map $ev_{S}:\cL(P_{2})\to\K^{n}$
satisfies the inequality in \rp{cond-rank}.
\item[5.] If yes, obtain $f_2$ with Newton polytope $P_2$ with coefficients from the matrix of the kernel of $ev_{S}$, stop.
If no, go back to Step 3 (or Step 1).
\end{itemize}

A variation of this algorithm is to loop over all irreducible polynomials in Step 1 until a toric complete intersection
is found. For $q=16$ this is still feasible. This is why in most our examples below $f_1$ does not look ``random".
In Step 3 we either run through all subsets of size $n$ or sample $10^5$ random subsets of size $n$, whichever is less, 
before we go back to Step~1.

We finish this section with a series of examples of toric complete intersection codes. 
Our first three examples demonstrate the construction in \rt{summands}, while the others
come from polygons in (1)--(3) of \rp{unmixed}. (We chose to only supply examples based on
polygons in the first row of \rf{fano}.)  The first example is written in full
detail, the reader may easily reconstruct details in the subsequent examples in a similar manner.

\begin{Ex} Our first example illustrates the construction of $P_1$, $P_2$ and $A$ in \rt{summands}.
We choose $R_1$ to be the vertical unit segment and $R_2$ a parallelogram ``around it", as in \rf{Th-2}. 
We put $A=P_1+R_1$ and $P_2=A+R_2$.
 \begin{figure}[h]
\includegraphics[scale=.7]{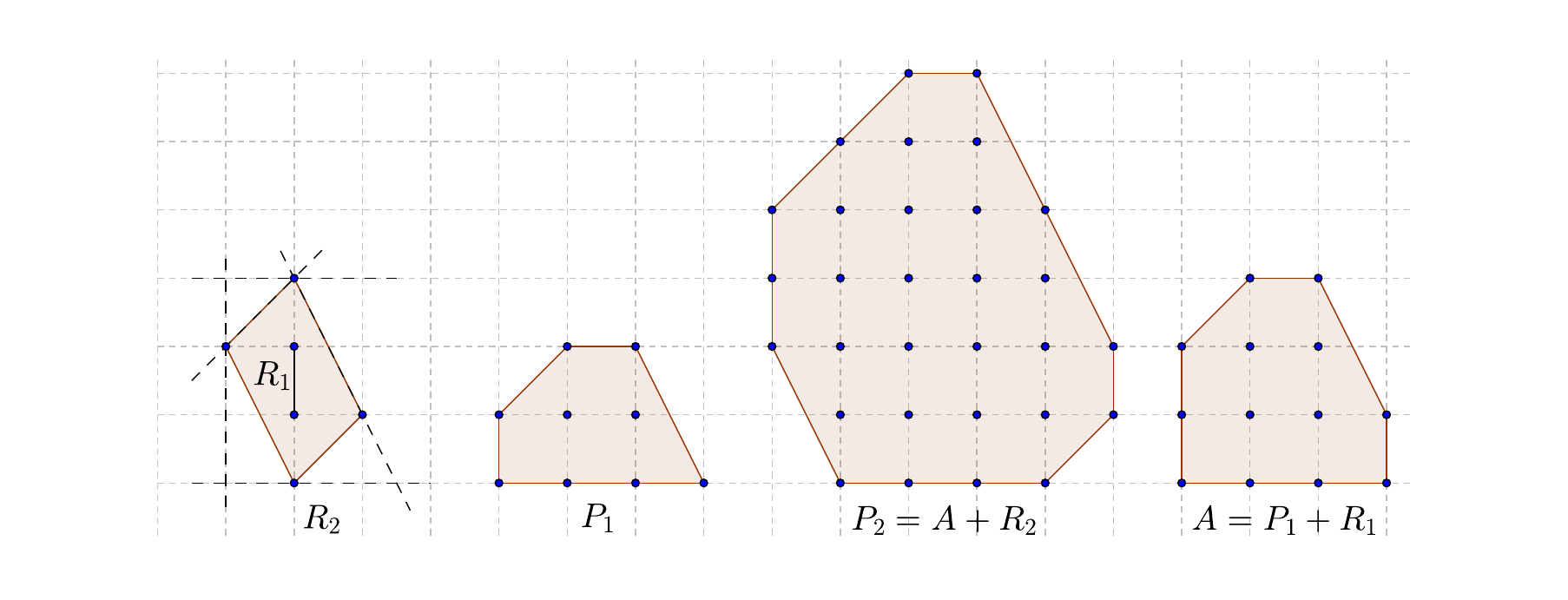}
\caption{Example of construction from \rt{summands}.}
\label{F:Th-2}
\end{figure}

Geometrically,  $l_{R_2}(v)= l_{R_1}(v)+1$, for all $v\in\Fan(P_1)$, means the following.  Draw lines parallel to the
sides of $P_1$ which are lattice distance one from $R_1$. We obtain, strictly speaking, a rational polygon
(presented by  dotted lines in \rf{Th-2}). Then the above condition means that $R_2$ is inscribed in this rational polygon.

Let $\K=\F_{16}$ with a primitive generator $t$. Consider the following system with Newton polygons $P_1$, $P_2$.

\begin{align*}
f_1 =&\, x^3 + x^2y^2 + t^7x^2y + t x^2 + xy^2 + x + y + 1=0,\\
f_2 =&\, t^7x^5y^2 + t^{11}x^5y + t^{11}x^4y^2 + t^9x^4y + t^7x^4 + t x^3y^6 + t^5x^3y + t^8x^3 + x^2y^6 + 
    t^5x^2y^3 + \\  &t^7x^2y^2 + t^5x^2y + t^{12}xy^5 + xy^4 + t^5xy^3 + xy^2 + t^{11}xy + t^{13}x + t y^4 + 
    t^4y^3 + t^9y^2=0.
\end{align*}

The solution set $S$ consists of $n=V(P_1,P_2)=22$ points in $(\F_{16}^*)^2$ and is a toric complete intersection. 
\begin{align*}
S=&\, \{ (1, t^{10}), (t, t^{10}), (t^3, 1), (t^3, t), (t^4, t^{10}), (t^5, t^7), (t^5, t^{11}), (t^6, t), 
(t^6, t^{13}), (t^7, t^2),\\ &\ \  (t^7, t^7), (t^8, t^2), (t^8, t^{12}), (t^9, t^{11}), (t^9, t^{13}), (t^{10}, 
t), (t^{10}, t^6), (t^{12}, t^7), (t^{12}, t^8), \\ & \ \ (t^{13}, t^{11}), (t^{13}, t^{14}), (t^{14}, t^4) \}.
\end{align*}
By \rp{dim}, $\dim(\cC_{S,A})=|A|_\Z-|R_1|_\Z=13-2=11$ which, as predicted by \rt{summands}, is exactly half the length of the code. By \rt{geom criterion}, $\cC_{S,A}$ is a quasi-self dual code. According to MAGMA, its parameters are $[22,11,10]$.
To find an equivalent self-dual code, first compute the vector $y$ of local residues:
\begin{align*}
y=&\,\Big(\dfrac{1}{J^{\mathbb{T}}_{f}(p_{1})},\dotsc,\dfrac{1}{J^{\mathbb{T}}_{f}(p_{22})}\Big)\\
=&\,\left( t^9, t^2, t, t^8, t^{11}, 1, t^{11}, 1, t, t^9, t^{10}, t^{14}, t^6, t^{10}, t^{12}, t^3, t^9, t^9, t^7, 
t^{11}, t^9, t^6 \right).
\end{align*}
This determines the vector $x$ such that $x_i^2=y_i$ for $1\leq i\leq 22$:
$$x=\left(t^{12}, t, t^8, t^4, t^{13}, 1, t^{13}, 1, t^8, t^{12}, t^5, t^7, t^3, t^5, t^6, t^9, t^{12}, t^{12}, 
t^{11}, t^{13}, t^{12}, t^3 \right).$$
Finally, the code $x\,\cC_{S,A}$ is a self-dual code with parameters $[22,11,10]$ over $\F_{16}$.

Next we look at some $y$-dual codes. Let $P^{(1)}$ denote the convex hull of the interior points of $P=P_1+P_2$.
Then we can decompose $P^{(1)}$ into Minkowski sum of two lattice polygons in several ways. They
are depicted in \rf{sums1}. (Of course, there is also $2A=P^{(1)}$, which we do not include in the figure.)

 \begin{figure}[h]
\includegraphics[scale=.7]{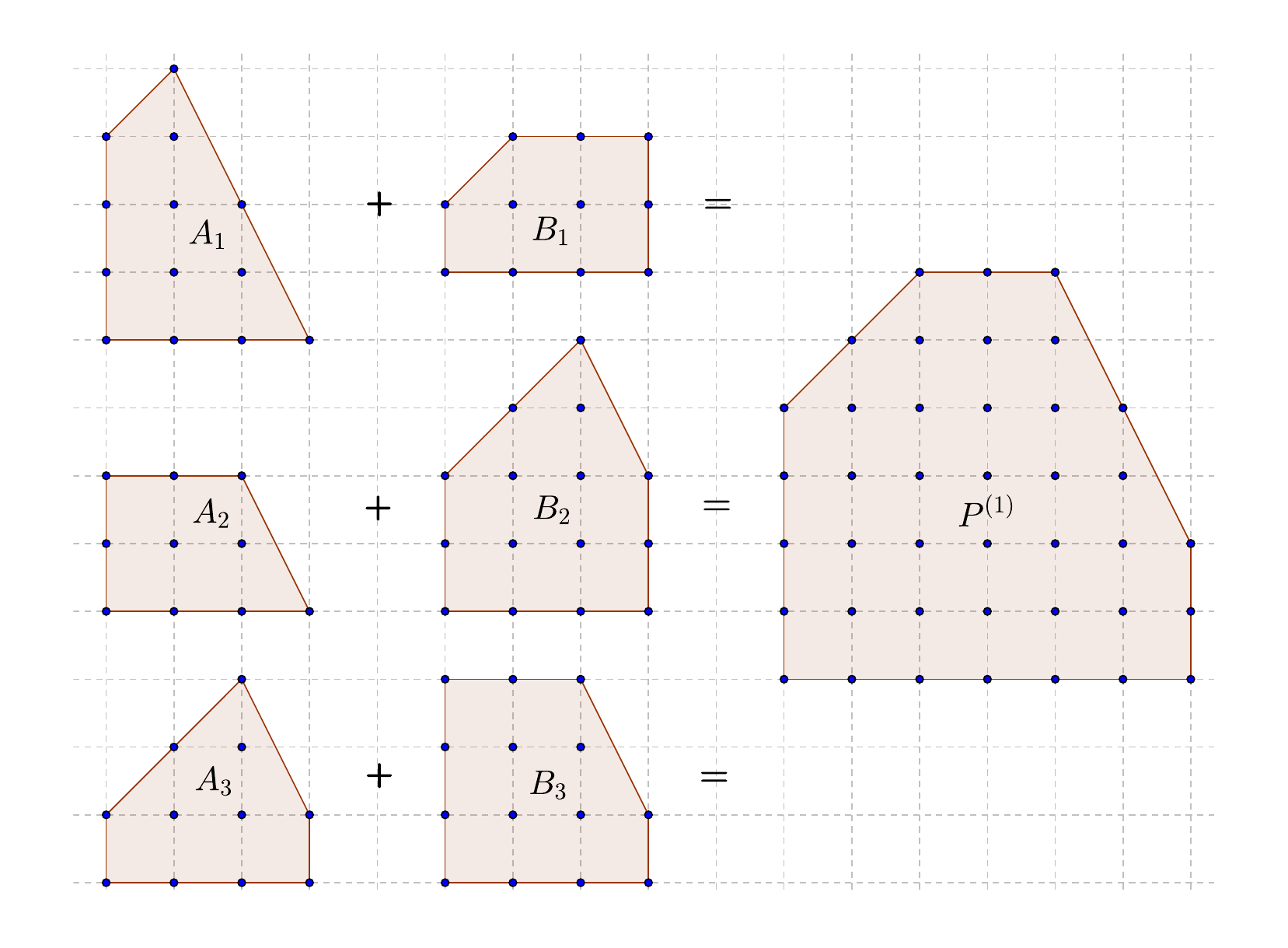}
\caption{Several Minkowski decompositions of $P^{(1)}$.}
\label{F:sums1}
\end{figure}
An easy application of \rp{dim} shows that the codes $\cC_{S,A_i}$ and $\cC_{S,B_i}$ have complementary dimensions.
Therefore, by \rt{y-dual toric}, they are $y$-dual codes. We list their parameters in a table below. 

 \bigskip
 
\renewcommand{\arraystretch}{1.5}
\begin{center}
    \begin{tabular}{ | c|p{2cm} | l | l | p{4.7cm} |} 
    \hline
    Polygons 
    & Parameters & Properties of Codes \\ \hline
$A_1$  
& $[22, 12, 9]$ &   \\
$B_1$  
& $[22,10, 11]$ & $y$-dual of $\cC_{S,A_1}$ \\ \hline
 $A_2$ 
 & $[22,9,12]$ &   \\ 
$B _2$ 
 & $[22,13,8]$ & $y$-dual of $\cC_{S,A_2}$ \\ \hline
  $A_3$ 
  & $[22,10,11]$ & \\
  $B _3$ 
& $[22,12,9]$ & $y$-dual of $\cC_{S,A_3}$   \\ \hline
  $A$ 
& $[22,11,10]$ & quasi-self-dual code  \\ \hline
    \end{tabular}
\end{center}

 \bigskip
 
\end{Ex}

\begin{Ex} Here is another illustration of \rt{summands}. This time $R_1$ is the 2-simplex and $R_2=P_1$ is 
a triangle containing $R_1$ in the interior (see \rf{Th-1}). As before, $A=P_1+R_1$ and $P_2=A+R_2$.

 \begin{figure}[h]
\includegraphics[scale=.7]{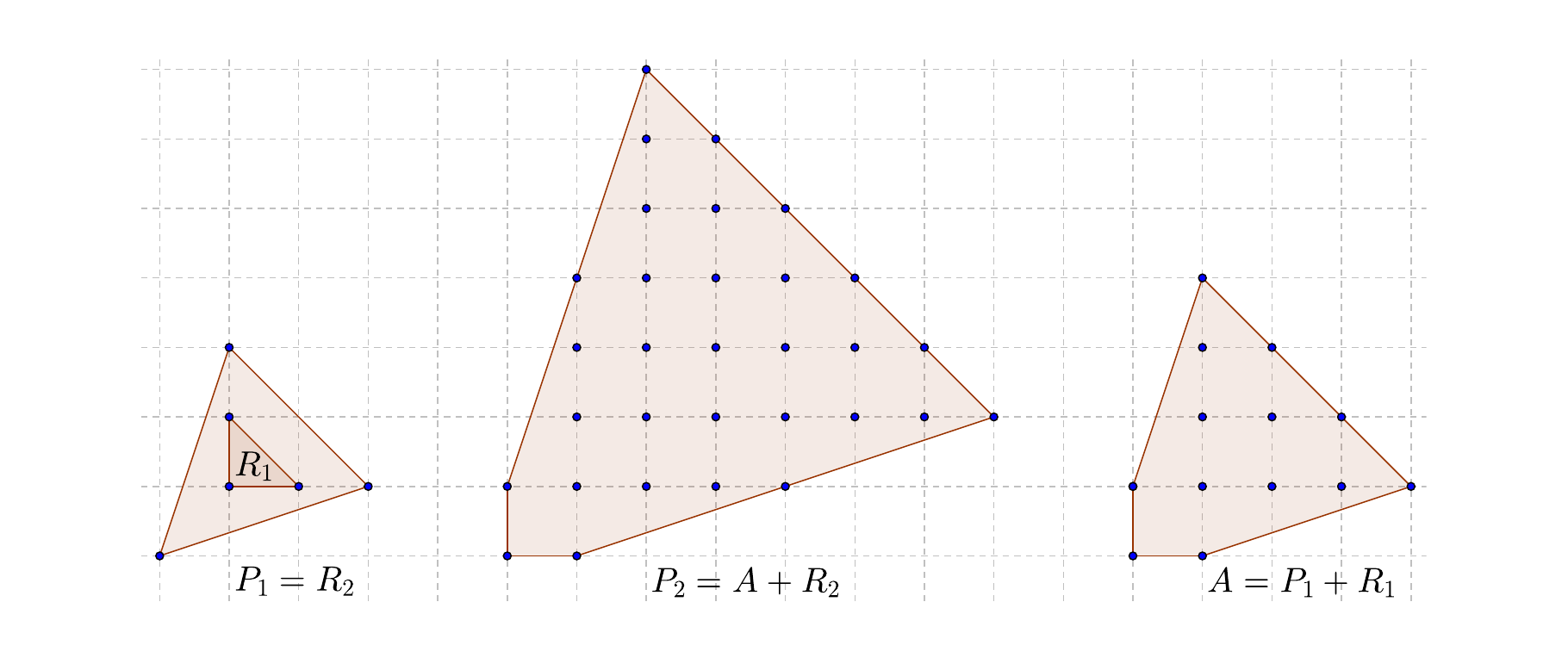}
\caption{Example of construction from \rt{summands}.}
\label{F:Th-1}
\end{figure}

Consider the following system over $\F_{16}$:
\begin{align*}
f_1= &\, x^3y + t^5x^2y^2 + x^2y + xy^3 + t^5xy^2 + xy + 1=0,\\
f_2= &\, t^5x^7y^2 + t^5x^5y^3 + t^5x^4y + x^3y^6 + t^5x^3y^2 + t^{10}x^3y + x^2y^7 + t^{10}x^2y + t^5x + t^5y + t^{10}=0.
\end{align*}

It defines a toric complete intersection $S$ of size $n=V(P_1,P_2)=20$. The quasi-self-dual code $\cC_{S,A}$
has parameters $[20,10, 8]$. As before, let  $P^{(1)}$ be  the convex hull of the interior points of $P$.
It is easy to check that $P_1$ is a Minkowski summand of $P^{(1)}$. We set $A_1=P_1$ and $B_1=P^{(1)}-A_1$.
The reader may wish to draw the polygons $A_1$, $B_1$, but we will omit it.
Below is the table of parameters of the corresponding $y$-dual codes.

 \bigskip
 
\renewcommand{\arraystretch}{1.5}
\begin{center}
    \begin{tabular}{ | c|p{2cm} | l | l | p{4.7cm} |} 
    \hline
    Polygons 
    & Parameters & Properties of Codes \\ \hline
$A_1$  
& $[20, 6, 12]$ &   \\
$B_1$  
& $[20,14, 4]$ & $y$-dual of $\cC_{S,A_1}$ \\ \hline
  $A$ 
& $[20,10, 8]$ & quasi-self-dual code  \\ \hline
    \end{tabular}
\end{center}

 \bigskip

\end{Ex}

\begin{Ex}
In our next example we consider rectangular boxes $P_1=[0,3]\times [0,2]$ and $P_2=[0,7]\times [0,4]$.
This is a particular case of \rt{summands}. First we choose a system  with these Newton polygons
which defines a toric complete intersection $S$ over  $\F_{16}$ of size $V(P_1,P_2)=26$:
\begin{align*}
f_1&\, = x^3y^2 + t^4x^3y + x^3 + t^5x^2y^2 + t^2x^2y + x^2 + t^{11}xy^2 + t xy + x + y^2 + y + 1 =0,\\
f_2&\, =x^7y^4 + t^{10}x^7y + t^2x^7 + t^{12}x^6y + t^8x^6 + t^{10}x^5y + t^3x^5 + t^{13}x^4y + t^{13}x^4 + t^9x^3y + t^{13}x^3 + \\  & 
t^{11}x^2y^3 +  t^8x^2y^2 + t^{12}x^2y + t^{14}x^2 + xy^4 + t^6xy^3 + t^3xy^2 + t^{12}x + t^6y^4 + t^9y^2 + t^{13}y + t^5  =  0.
\end{align*}

Let  $P^{(1)}$ be the convex hull  of the interior lattice points  of $P$, shifted to the origin, i.e. $P^{(1)}=[0,8]\times, [0,4]$.
Then, a shift of $A=[0,4]\times [0,2]$ defines a quasi-self-dual code. 
We also try different subsets $A$ and $B$ such that $A+B=P^{(1)}$. Note that $A+B=P^{(1)}$
does not guarantee that $\cC_{S,A}$ and $\cC_{S,B}$ are $y$-dual since their dimensions might
not be complementary. For example, if $A=[0,3]\times [0,4]$ and $B=[0,5]\times \{0\}$ then $A+B=P^{(1)}$.
However, $\dim\cC_{S,A}=20-3=17$ and $\dim\cC_{S,B}=6$, by \rp{dim}.
In the table below all the codes have best known parameters as confirmed in \cite{mint}. 
\bigskip
 
\renewcommand{\arraystretch}{1.5}
\begin{center}
    \begin{tabular}{ | p{3.5cm} | l | l | p{4.7cm} |} 
        \hline
    Polytopes 
    & Parameters & Properties of Codes \\ \hline
 (a) $A =[0,3]\times[0,3]$ 
 & $[26,14,11]$ & 
 \\ 
  \hspace{.45cm}  $B =[0,5]\times[0,1]$  
  & $[26,12,13]$ & 
  $y$-dual of $\cC_{S,A}$ \\ \hline
 (b) $A=[0,4]\times[0,2]$ 
 & $[26,13,12]$ & quasi-self-dual code  \\ \hline
 (c) $A=[0,2]\times[0,2]$ 
 & $[26,9,16]$ &   \\ 
 \hspace{.45cm}  $B =[0,6]\times[0,2]$ 
 & $[26,17,8]$ & $y$-dual of $\cC_{S,A}$  \\ \hline
 (d) $A=[0,3]\times[0,1]$ 
 & $[26,8,17]$ &   \\ 
 \hspace{.45cm}  $B =[0,2]\times[0,4]$ 
 & $[26,18,7]$ &$y$-dual of $\cC_{S,A}$  \\ \hline
     \end{tabular}
\end{center}
 \bigskip

\end{Ex}

The next series of examples uses polygons classified in \rp{unmixed}. We construct dual and quasi-self-dual
toric complete intersection codes for cases (1), (2), and the first four polygons
in case (3).

\begin{Ex} In this example $P_1=3\,\triangle$, $P_2=6\,\triangle$, where $\triangle$ is the standard $2$-simplex.
Consider the following system over $\F_{16}$:
\begin{align*}
f_1 =&\, x^3 + x^2y + x^2 + t^3xy^2 + xy + x + y^3 + t y^2 + y + 1= 0, \\
f_2 =&\, t^{11}x^6 + t^3x^5y + t^9x^5 + t x^4y^2 + t x^4y + t^{13}x^4 + t^6x^3y + t^9x^3 +  t^4x^2y^2 + t^3x^2y 
    + \\ & t^8x^2 + t^2xy^2 + t^4xy + y^6 + t^8y + t^6  =  0.
\end{align*}
The solution set $S$ is a toric complete intersection consisting of $V(P_1,P_2)=3\cdot 6=18$ points. 
Below is a table of codes for different choices of $A$ and $B$ whose sum lies in $P^{\circ}$
up to a lattice translation.


 \bigskip
 
\renewcommand{\arraystretch}{1.5}
\begin{center}
    \begin{tabular}{ | p{2.5cm} | l | l | p{4.7cm} |} 
    \hline
    Polytopes 
    & Parameters & Properties of Codes \\ \hline
(a) $A =\triangle$  
& $[18,3,15]$ &   \\
  \hspace{.45cm}  $B =5\,\triangle$  
& $[18,15, 3]$ & $y$-dual of $\cC_{S,A}$ \\ \hline
 (b) $A =2\,\triangle$ 
 & $[18,6,12]$ &   \\ 
  \hspace{.45cm}   $B =4\,\triangle$ 
 & $[18,12,6]$ & $y$-dual of $\cC_{S,A}$ \\ \hline
  (c) $A=3\,\triangle$ 
& $[18,9,9]$ & quasi-self-dual code  \\ \hline
    \end{tabular}
\end{center}

 \bigskip

\end{Ex}

\begin{Ex} Now let $P_1=3\,\square$, $P_2=5\,\square$, where $\square=[0,1]\times[0,1]$ is the unit square.
Let $\K=\F_{16}$ as before, and consider the following system:

\begin{align*}
f_1=&\, x^3y^3 + x^3y + x^3 + x^2y^2 + x^2y + x^2 + xy^3 + xy^2 + xy + x + y^3 + y^2 + y + 1=0.\\
f_2=&\, t x^5y^5 + x^5y^2 + x^5y + t x^5 + x^4y^5 + t^4x^4y + x^4 + x^3y^2 + x^3y + x^2y^3 + x^2y +\\ &  t xy^4 + x + t y^5 + y^4 + y + 1=0.
\end{align*}
The solution set $S$ consists of $V(P_1,P_2)=2\cdot 3\cdot 5=30$ points and is a toric complete
intersection. Here is a table of codes for different choices of $A$ and $B$.

 \bigskip
 
\renewcommand{\arraystretch}{1.5}
\begin{center}
    \begin{tabular}{ | p{3.5cm} | l | l | p{4.7cm} |} 
        \hline
    Polytopes 
    & Parameters & Properties of Codes \\ \hline
 (a) $A =\square$ 
 & $[30,4,24]$ & 
 \\ 
  \hspace{.45cm}  $B =5\,\square$  
  & $[30,26,3]$ & 
  $y$-dual of $\cC_{S,A}$ \\ \hline
 (b) $A=3\,\square$ 
 & $[30,15,12]$ & quasi-self-dual code  \\ \hline
 (c) $A=2\,\square$ 
 & $[30,9,18]$ &   \\ 
 \hspace{.45cm}  $B =4\,\square$ 
 & $[30,21,6]$ & $y$-dual of $\cC_{S,A}$  \\ \hline
 (d) $A=[0,4]\times[0,2]$ 
 & $[30,15,12]$ &   \\ 
 \hspace{.45cm}  $B =[0,2]\times[0,4]$ 
 & $[30,15,12]$ &$y$-dual of $\cC_{S,A}$  \\ \hline
  (e) $A=[0,3]\times[0,2]$  
  & $[30,12,15]$ & 
  \\
   \hspace{.45cm}  $B =[0,3]\times[0,4]$ 
   & $[30,18,9]$ & $y$-dual of $\cC_{S,A}$\\ \hline
    \end{tabular}
\end{center}
 \bigskip

\end{Ex}

\begin{Ex}
Let $P_1=2Q$ and $P_2=4Q$ where $Q$ is the convex hull of $\{0,2e_1, e_2\}$.
Consider the following system over $\F_{16}$.

\begin{align*}
f_1=&\, x^4 + x^3 + t x^2y + x^2 + xy + x + y^2 + t y + 1=0,\\
f_2=&\, t^7x^8 + t^{13}x^7 + t^{12}x^6 + t^{14}x^5y + t^2x^5 + t^9x^4y + t^9x^4 + t^5x^3y + t^{11}x^3 +\\ &
 t^8x^2y + t^8x^2 + t^{11}xy + x + y^4 + t^{12}y + t^{13}=0.
\end{align*}

Its solution set $S$ is a toric complete intersection of size $16$. For $A=Q$, $2Q$, and $3Q$
the corresponding codes $\cC_{S,A}$ have dimensions $4$, $8$, and $14$, respectively, by \rp{dim}.
We have the following table of codes and their parameters.

 \bigskip
 
\renewcommand{\arraystretch}{1.5}
\begin{center}
    \begin{tabular}{ | c|p{2cm} | l | l | p{4.7cm} |} 
    \hline
    Polygons 
    & Parameters & Properties of Codes \\ \hline
(a) $A =Q$  
& $[16, 4, 12]$ &   \\
  \hspace{.45cm} $B=3Q$  
& $[16,12, 4]$ & $y$-dual of $\cC_{S,Q}$ \\ \hline
(b) $A =2Q$ 
& $[16,8, 8]$ & quasi-self-dual code  \\ \hline
    \end{tabular}
\end{center}

\bigskip

\end{Ex}

\begin{Ex}
Let $P_1=2Q_1$ and $P_2=3Q_1$ where $Q_1$ is the first Fano polygon as in \rf{fano} and consider the following system
with these Newton polygons.

\begin{align*}
f_1=&\, x^4y^4 + t^5x^3y^2 + t^{10}x^2y^3 + x^2y^2 + x^2y + x^2 + xy^2 + xy + y^2=0,\\
f_2=&\, x^6y^6 + x^4y^4 + x^4y^3 + x^3y^4 + x^3y^3 + t^5x^3y^2 + t^5x^3y + t^5x^3 + \\ & t^{10}x^2y^3 + t^{10}xy^3 + t^{10}y^3=0.
\end{align*}

Its solutions set $S$ is a toric complete intersection of size $n=2V(Q_1)\cdot 6=18$. By \rp{unmixed},
the code $\cC_{S,A}$ with $A=2Q_1$ is quasi-self dual. Now, notice that $Q_1$ and $3Q_1$ satisfy
$Q+3Q\subseteq P^\circ$. In fact,  the corresponding codes have complementary dimensions.
Indeed, $\dim(\cC_{S,Q_1})=|Q_1|_\Z=4$, clearly.  
As for $\dim(\cC_{S,3Q_1})$, \rp{dim} is not applicable since $V(A,P_2)<V(P_1,P_2)$ fails. 
But it's clear here that the kernel of the evaluation map $ev_{S}:\cL(3Q_1)\to\F_{16}^{18}$ has one more basis element,
namely, $P_2$ itself. Therefore, $\dim(\cC_{S,3Q_1})=|3Q_1|_\Z-|Q_1|_\Z-1=19-4-1=14$. 
This justifies that $\cC_{S,Q_1}$ and $\cC_{S,3Q_1}$ are $y$-dual. We record the corresponding parameters below.
 
 \bigskip
 
\renewcommand{\arraystretch}{1.5}
\begin{center}
    \begin{tabular}{ | c|p{2cm} | l | l | p{4.7cm} |} 
    \hline
    Polygons 
    & Parameters & Properties of Codes \\ \hline
(a) $A = Q_1$  
& $[18, 4, 13]$ &   \\
  \hspace{.45cm} $B=3Q_1$  
& $[18,14, 4]$ & $y$-dual of $\cC_{S,Q_1}$ \\ \hline
(b) $A =2Q_1$ 
& $[18,9, 8]$ & quasi-self-dual code  \\ \hline
    \end{tabular}
\end{center}

\bigskip

Similarly, we obtain toric complete intersection for $Q_2$, $Q_3$, and $Q_4$ (the polygons in the first row of \rf{fano}).
For $P_1=2Q_2$ and $P_2=3Q_2$ we take 
\begin{align*}
f_1=&\, x^4 + x^3y^2 + x^3y + x^3 + x^2y^4 + t^8x^2y^3 + t^3x^2y^2 + x^2y + x^2 + xy^2 + xy + x + 1=0,\\
f_2=&\, t^8x^6 + t^{13}x^5y^2 + t^{13}x^5y + x^5 + x^4y^2 + t^2x^4y + t^7x^4 + x^3y^6 + t^{13}x^3y^3 + t^{14}x^3y^2 + \\ & t^{13}x^3y + 
    t^2x^3 + x^2y^2 + t^2x^2y + t^7x^2 + t^{13}xy^2 + t^{13}xy + x + t^8=0.
\end{align*}

For $P_1=2Q_3$ and $P_2=3Q_3$ we take 
\begin{align*}
f_1=&\, x^4y^2 + t^{11}x^3y^3 + x^3y^2 + x^3y + x^2y^4 + t^9x^2y^3 + x^2y^2 + x^2y + x^2 + t^{11}xy^3 + xy^2 + xy + y^2=0,\\
f_2=&\, t^5x^6y^3 + t^4x^5y^3 + t^2x^5y^2 + t^6x^4y^4 + x^4y^3 + t^{12}x^4y^2 + t^7x^4y + x^3y^6 + t^{14}x^3y^3 + \\ &
    t^7x^3y^2 + tx^3y + t^{13}x^3 + t^6x^2y^4 + x^2y^3 + t^{12}x^2y^2 + t^7x^2y + t^4xy^3 + t^2xy^2 + t^5y^3=0.
\end{align*}

Finally, for $P_1=2Q_4$ and $P_2=3Q_4$ we take
\begin{align*}
f_1=&\, x^4y^4 + t^8x^3y^3 + t^{13}x^3y^2 + t^{13}x^2y^3 + t^2x^2y^2 + x^2y + x^2 + xy^2 + xy + x + y^2 + y + 1=0,\\
f_2=&\, x^6y^6 + t^7x^4y^2 + t^{14}x^3y^3 + t^2x^3y^2 + t^7x^3y + tx^3 + t^7x^2y^4 + t^2x^2y^3 + t^8x^2y^2 + 
    t^2x^2y\\ & + t^3x^2 + t^7xy^3 + t^2xy^2 + t^{14}xy + t^4x + ty^3 + t^3y^2 + t^4y + t^{10}=0.
\end{align*}

The corresponding codes  happen to have the same parameters for each $i=2,3,4$ and are listed below.

 \bigskip
 
\renewcommand{\arraystretch}{1.5}
\begin{center}
    \begin{tabular}{ | c|p{2cm} | l | l | p{4.7cm} |} 
    \hline
    Polygons 
    & Parameters & Properties of Codes \\ \hline
(a) $A =Q_i$  
& $[24, 5, 16]$ &   \\
  \hspace{.45cm} $B=3Q_i$  
& $[24,19, 4]$ & $y$-dual of $\cC_{S,Q_i}$ \\ \hline
(b) $A =2Q_i$ 
& $[24,12, 8]$ & quasi-self-dual code  \\ \hline
    \end{tabular}
\end{center}

\bigskip

\end{Ex}

\begin{Ex} Our final example does not use the geometric construction of \rt{summands}. The polygons are depicted in \rf{general}.
 \begin{figure}[h]
\includegraphics[scale=.72]{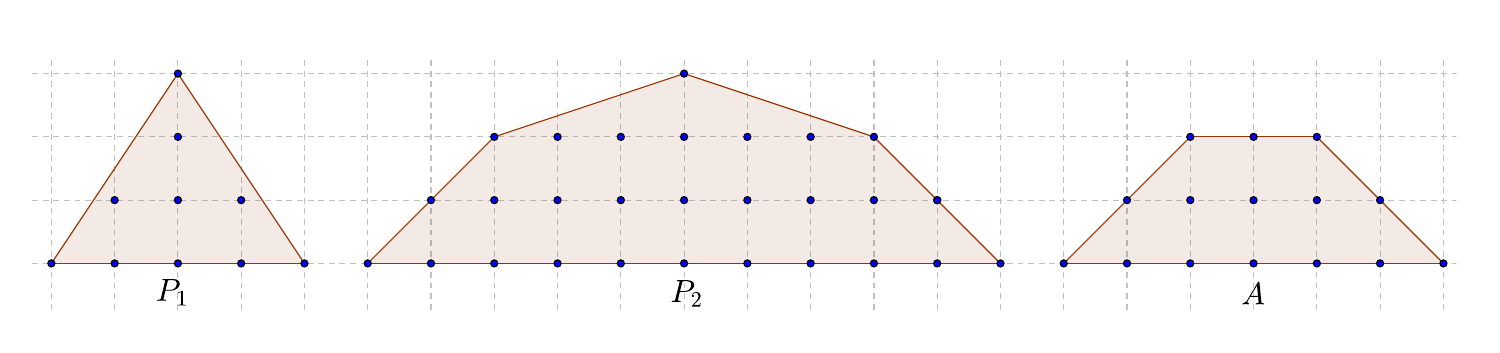}
\caption{Here $P_1$ is not a Minkowski summand of $A$.}
\label{F:general}
\end{figure}
Clearly, $P_1$ is not a Minkowski summand of $A$,
but one can check that the equality of the support functions in \rr{general} still holds. The toric complete intersection $S$ is defined over $\F_{16}$
by
\begin{align*}
f_1=&\, x^4 + x^2y^3 + x^2 + x + 1=0,\\
f_2=&\, t^8x^{10} + tx^9y + t^3x^9 + t^9x^8y^2 + t^{10}x^8y + t^{12}x^8 + t^4x^7y^2 + t^{14}x^7y + t^4x^7 + t^7x^6y^2 +\\ & t^6x^6y + tx^6 + t^9x^5y^3 +  tx^5y + t^4x^5 + t^7x^4y^2 + t^7x^4y + t^6x^4 + t^{10}x^3y^2 + x^3y + t^6x^3 + x^2y^2\\ & + x^2y + t^6x^2 + xy + t^{10}x + 1=0.
\end{align*}
The corresponding quasi-self-dual code $\cC_{S,A}$ has parameters $[30,15,12]$.
\end{Ex}

\end{document}